\documentclass[10pt,a4paper]{article}

\usepackage{latexsym, amssymb, amsfonts, amsmath}
\usepackage{graphicx}
\usepackage{frege}

\usepackage[numbers]{natbib}
\usepackage[displaymath, mathlines]{lineno}

\newcommand{\mh}{\noindent}
\newcommand{\cao}{\c{c}\~{a}o}

\newcommand{\negr}[1]{\boldsymbol{#1}}

\begin{document}

\title{\bf Frege on the reference of sentences} 
\author{Ab\'ilio Rodrigues \\ 
[4mm]
Federal University of Minas Gerais \\  
abilio.rodrigues@gmail.com}

\date{}

\maketitle

\begin{abstract}
\mh The aim of this paper is to show that Frege's argument which concluded that the reference of a sentence is its truth-value, presented in \textit{On Sense and Reference} (1892), can be reconstructed taking into account the problems of the notion of conceptual content presented in the \textit{Begriffsschrift} (1879) and also other passages from a letter to Russell (1902) and the posthumous \textit{Logic in Mathematics} (1914). 
Once the `hybrid' notion of conceptual content 
was rejected as the semantic value of the expressions of the formal language 
designed to carry out the logicist project,  there was no alternative between truth-values and thoughts.
I claim that the reconstructed argument is perfectly sound and convincing.
\end{abstract}



\section*{Introduction}

The thesis that the reference of a sentence is its truth-value is a central point 
of Frege's work in Logic. It was criticized for the first time by Russell in a letter to Frege in 1902. 
In this letter, Russell objected that he could not 

\begin{quote}
believe that the true or the false is the reference of a proposition in the same sense as, e.g., a certain person is the reference\footnote{In \cite{frege.corresp} we read `meaning' instead of `reference'. The discussion about the translation of `Bedeutung' is well-known. I have changed here all the occurrences of `meaning' rendering `Bedeutung' by `reference', in quotations from \citep{transl, frege.posthu, frege.corresp}.} of the name Julius Caesar \citep[pp. 150-151]{russell}.     
\end{quote}

\mh 
Russell's view is somewhat similar to both Dummett's  \citep{dummett73} and Chateaubriand's 
\citep{chat} views. 
According to Dummett, we are justified in taking the relation between a proper name and its bearer as the prototype for the relation between a sentence and its reference, and so he argues that Frege should not have ascribed reference to sentences \citep[p. 181]{dummett73}.   
At this point, Dummett is doubly mistaken: the prototype is not the relation between a name and its bearer  and, of course, Frege had to ascribe reference (i.e., semantic value) to the sentences of his formal language. Some pages later, we read:

\begin{quote}

The identification of truth-values as referents of sentences, taken together with the thesis that the truth-values are objects, led to a great simplification of Frege's ontology, at the price of a highly implausible analysis of language. (...) It is tragic that a thinker who achieved the first really penetrating analysis of the structure of our language should have found himself driven to such absurdities (...) the assimilation of sentences to proper names did have a fatal effect upon Frege's theory of meaning. 
It is just that Frege's earlier departures from the forms of natural language (...) were founded upon deep insights into the workings of language; whereas this ludicrous deviation is prompted by no necessity, but is a gratuitous blunder \citep[p. 184]{dummett73}.
\end{quote}
\mh
As it is well-known, Dummett reads Frege from the viewpoint of his (Dummett's) 
interest in a theory of meaning.  
Frege, however,  
was interested in devising a formal language 
and an account of logical consequence to carry out his logicist project. 
There is no `fatal effect upon Frege's theory of meaning', as Dummett claims,  because Frege  
did not have a theory of meaning, nor was he interested in an analysis of natural language.
Notice, likewise, that it 
is entirely possible to understand Frege's claim that truth-values are objects as pragmatically 
motivated, similar to our practice of writing on the blackboard things like $ I(p) = F $, 
$ I(\neg p) = T $, and $I(a) = Aristotle$ when teaching logic -- 
of course, in writing such things, we are taking sentences as names and truth-values as objects.
We will see that that the `identification of truth-values as referents of sentences' is anything but a gratuitous blunder. 
Instead, it is fair to say (paraphrasing Dummett)  that it is a result of necessity,
founded upon a deep insight  into the workings of a formal language for Mathematics.

\

Chateaubriand also criticizes Frege's choice, in my view, committing the same error as Dummett:

\begin{quote}
[I]f one looks for candidates, one will hit upon states of affairs, facts, or something like that. 
The similarity between a definite description, as a way of presenting an object, 
and a sentence, as a way of presenting a state of affairs, is so striking and obvious (...) 
that it just can't be missed \citep[p. 76]{chat}.  
\end{quote}

\mh All these criticisms are mistaken.  
The most important point  is that Frege just {realized} -- or discovered --  that:
(i)  an extensional logic is very suitable as an account of logical consequence for Mathematics, 
and  (ii) in such an extensional logic,  
the extension of a sentence is its truth-value. 
The fact that Frege informally explained the notion of reference starting from the puzzle of 
identity and based on examples of natural language, where the relation between a singular 
term and its 
reference is that of a name and its bearer, should not be a problem. 
Frege's strategy may be understood as 
an informal and pre-theoretical elucidation, 
in the sense explained by Weiner in  \citep{weiner}. The relation between a 
sentence and a truth-value is not like the relation between a name and its bearer. 
Rather, it is a relation between 
a linguistic expression and its semantic value within a formal system. The relation between a singular 
term and the corresponding object is also a relation between an expression and its semantic value, 
although, in this case, it also happens to be a relation between a name and its bearer.

Among his published works, the one in which Frege argues in defense of the thesis that the reference 
of a sentence is its truth-value is the paper \textit{On Sense and Reference} \cite{SR} (from now on 
\textit{SR}), published in 1892, one year before the first volume of \textit{The Basic Laws of Arithmetic} 
\citep{bla} (from now on \textit{BLA}). The latter would be the main work of Frege's academic career, were 
it not for the inconsistency of the fateful Basic Law V. The purpose of the papers published in 1891 and 
1892, \textit{Function and Concept} \citep{FC} and \textit{SR}, was to present some adjustment in the 
formal system to be used in \textit{BLA}. In particular, the primary aim of \textit{SR} was to establish 
that the reference (or better, the extension) of a sentence is its truth-value. This thesis, together with 
the sense/reference distinction for singular terms and Frege's theory of extensions (that ended up being 
inconsistent), yields an extensional logic for \textit{BLA}. 
On the other hand, it is  true that in 1892 Frege's argument for truth-values as references had some 
problems, and it may seem that he was not wholly convinced. 
	
In fact, the above narrative has already been established as the standard interpretation of \textit{SR}, 
but, in my view, there are some points of Frege's line of reasoning that are not yet fully explored. 
The aim of this paper is to show that Frege's argument may be reconstructed taking into account 
the failure of the notion of conceptual content, which in the \textit{Begriffsschrift} \citep{begr} 
(from now on \textit{BS}) played the role of semantic value of the expressions of his formal system, 
and also other passages where Frege presents the argument with a small but important difference: 
a letter to Russell from 1902 \citep{letter} and the text \textit{Logic in Mathematics} \citep{LM}, 
published posthumously and dated 1914 by the editors of \cite{frege.posthu}. 
My claim here is that the reconstructed argument is perfectly sound and convincing.
	
\
	
The remainder of this text is structured as follows. In Section \ref{sem.values}, I recall some notions 
useful  for the discussion to be carried out here, namely, the notion of semantic value in a formal 
system, and the distinction between intension and extension. In Section \ref{conc.cont}, 
I discuss the problems of the notion of conceptual content presented by Frege in 1879 in \textit{BS}. 
In Section \ref{recon}, Frege's argument in defense of truth-values as references is analyzed and 
reconstructed, based not only on \textit{SR} but also on the letter to Russell  \citep{letter} 
and on the posthumous \textit{Logic in Mathematics} \citep{LM}. 

\section{Semantic value, extension, intension } \label{sem.values}
 
In a formal system, the semantic value of a linguistic expression is an entity associated with that linguistic expression. Unless we are talking about a  language, the semantic value is something 
non-linguistic.  
In compositional semantics, the semantic value of a complex expression depends functionally on the semantic values of their constitutive parts and on the way they are combined. Let \textit{v} be the semantic value of an expression \textit{A}. When \textit{A} is part of a more complex expression $(...A...)$ the semantic value of $(...A...)$ depends on \textit{v} and on the structure of $(...A...)$. Moreover, if $(...A...)$ is a sentence, 
the truth-value of $(...A...)$ also 
depends on \textit{v} and on the structure of $(...A...)$.

The semantic value associated with an expression \textit{A} may be the extension or the intension of \textit{A}. Roughly speaking, the intension of an expression \textit{A} is the meaning of \textit{A}, and the extension is what is denoted by \textit{A}.\footnote{More detailed explanations of the intension/extension distinction can be found  in \citep{fitt} and~\citep[pp. 3-14]{kirkham}.} 
The singular terms `the author of Nicomachean Ethics' and `the tutor of Alexander the Great' have the same extension, namely, Aristotle, but their meanings are not the same. 
Concerning predicates, this distinction is  illustrated by the typical example of the predicates `\textit{x} is a human being' and `\textit{x} is a featherless biped'. 
Both have the same extension because the set associated to each one is the same, but the intensions, that is, the meanings of the predicates, are not the same. The intension of a sentence is usually taken to be the proposition expressed, and the extension is its truth-value --  precisely  Frege's thesis discussed here.

 A logic $\mathcal{L}$ with  a language \textit{L} is  intensional  when $\mathcal{L}$ is  concerned not only with the extensions  but also with the intensions  of the expressions of \textit{L}. 
 Otherwise, if $\mathcal{L}$ is concerned only with the extensions of the expressions of \textit{L}, 
we say that  $\mathcal{L}$ is extensional. 
In other words, a logic is extensional when  it is not concerned with  the meanings of its expressions, but only with the entities referred  by them, no matter how these entities are picked up.
Analogously, a logical operator $ \phi $ is said to be extensional if its behaviour depends 
only on the extensions of the expressions in its scope. 
Classical first-order logic is {extensional}.
 It has been designed to formalize mathematical reasoning, and it is very well suited to this task.
 Classical first-order logic is not a theory of meaning, nor it has been primarily conceived to formalize  argumentative contexts of everyday reasoning.
 
For singular terms and sentences, the intension/extension distinction corresponds to the Fregean 
sense/reference distinction.\footnote{It is worth noting that there is a point that may 
cause some confusion to the reader of \textit{SR}. 
Although Frege distinguishes two semantic aspects associated with a linguistic expression, the reference and the sense (i.e., the extension and the intension) the word `reference' sometimes does not mean the extension.   
It happens when Frege says that in indirect speech the reference of a sentence is the thought.
 In more precise terminology, we would say that `the reference of a sentence is its truth-value and the sense is the thought' means, `the extension of a sentence is its truth-value and the intension is the thought', and  `in indirect speech, the reference of a sentence is the thought' means `the semantic value of a sentence in an intensional context is the thought'.} In a posthumously published paper \citep{CSB}, Frege says that the reference of a predicate is a concept, and that the sense is its mode of presentation, but he doesn't 
explain what would be a mode of presentation of a concept. In a letter to Husserl \citep{letter.h}, written in 1891, however, Frege makes it clear that he conceived of concepts extensionally.

\section{Conceptual content in the \textit{Begriffsschrift} (\textit{BS})} \label{conc.cont}

At the time of the \textit{Begriffsschrift} \citep{begr} (1879), Frege had not established  the distinction between sense and reference yet. In that work, he talks about conceptual content and judgeable content. The latter qualifies the content of sentences, that is, a content that is either true or false, and so it is a particular case of the former. 
I will talk here only about the more general notion, conceptual content, that applies to both singular terms  and sentences. 

\subsection{Conceptual contents of sentences}

Conceptual contents of sentences  are introduced in the 
\S3 of \textit{BS} and explained in terms of inferential role:
 
\begin{quote}
the contents of two judgments can differ in two ways: either the consequences derivable from the first, when it is combined with certain other judgments, always follow also from the second, when it is combined with these same judgments, [and conversely] or this is not the case. The two propositions `The Greeks defeated the Persians at Plataea' and `The Persians were defeated by the Greeks at Plataea' differ in the first way. Even if one can detect a slight difference in meaning, the agreement outweights it.  
Now I call that part of the content that is the \textit{same} in both the \textit{conceptual content}.
Since it \textit{alone} is of significance for our ideography [Begriffsschrift], we need not introduce 
any distinction between propositions having the same conceptual content. (...) 
   
\mh
[I]n a judgment I consider only that which influences its possible consequences. Everything necessary for a correct inference is expressed in full. \citep[\S3]{begr}.
\end{quote}

The conceptual content of a judgment (or a sentence) is what is relevant for inferences. 
Thus, we begin by saying that two sentences \textit{A} and \textit{B} have the same conceptual content if and only if they 
are intersubstitutable preserving correctness of inference. 
However, for Frege, correctness of inference was not exactly what we understand today by logical consequence. 
The former is not a sufficient nor a necessary condition for the latter. 
If all true arithmetical propositions are logical truths, as Frege held, then they will have the same conceptual content, but  Frege certainly would not agree with this. 
Moreover, the example mentioned by Frege of a relation  and its inverse  does not qualify as logical equivalence, 
since in first-order logic, $Rab \leftrightarrow R^{-1}ba$ is not a logical truth. 
So, logical equivalence cannot be a criterion of identity for conceptual contents.  

	Let us consider for the sake of the argument that  $\Gamma \vdash_F A $  means that through  one or more `Fregean correct inferences' the sentence $A$ may be obtained from the set of sentences $ \Gamma $, 
that is, the symbol $ \vdash_F $ here means `Fregean logical consequence'. 
So, from the passage quoted above, we get the following criterion for sameness of conceptual content of sentences:

\begin{description}
\item[] 	(1) $A$ and $B$ have the same conceptual content if and only if  for any $ \Gamma $ and $ C $: $ \Gamma, A \vdash_F C $ iff  $ \Gamma, B \vdash_F C  $. 
\end{description}

\mh With respect to sentences, and according to \S3 of \textit{BS}, the identity of content between 
\textit{A} and \textit{B} works analogously to  the replacement property,  according to which equivalent formulas are 
logically indistinguishable:
\begin{description}
\item[] If $A \dashv \vdash B$, then $ C(A/p) \dashv \vdash  C(B/p)$,
where $C(A/p)$ results from substituting one or more occurrences of \textit{p} by \textit{A} in
\textit{C}.
\end{description}
\mh  
Indeed, if we assume that the relation $ \vdash_F $ is Tarskian\footnote{A given 
relation of consequence $ \vdash $ is Tarskian, or standard, when  reflexivity, monotonicity, and transitivity hold  for $ \vdash $.}, the criterion (1) above 
is equivalent to 

\begin{description}
\item 	(2) $A$ and $B$ have the same conceptual content if and only if $  A \vdash_F B $ and  $ B \vdash_F A $.
\end{description}

\subsection{Conceptual content of singular terms}

\mh 
The natural way of extending this idea to singular terms is by saying that two singular terms \textit{a}
and \textit{b} have the same conceptual content when they are intersubstitutable preserving correctness of inferences, and so, 

	\begin{description}
\item (3) Two singular terms \textit{a} and \textit{b} have the same conceptual content if and only~if 	
	$ (\dots a \dots) \vdash_F (\dots b \dots) $ and  $ (\dots b \dots) \vdash_F (\dots a \dots) $.
	\end{description}

\mh
We will see, however, that the notion of conceptual content so defined is not compatible with Frege's notion of conceptual content for singular terms.

\

In \S8 of \textit{BS} Frege presents the sign of identity of content $ \equiv $. 
This sign is to be understood metalinguistically, that is, $ A \equiv B$ expresses a relation between the signs \textit{A} and \textit{B}, namely, that they have the same content. He takes an example from geometry in which a point is determined in two different ways and introduces a distinction that is virtually the same as the subsequent distinction between sense and reference for singular terms.
	 
\begin{quote}
To each of these ways of determining the point there corresponds a particular name. Hence, the need for a sign for identity of content rests upon the following consideration: the same content can be completely determined in different ways; but that in a particular case \textit{two ways of determining it} 
really yield the same result is the content of a \textit{judgment}. Before this judgment can be made, two distinct names corresponding to the ways of determining the content, must be assigned to what these ways determine. Before this judgment can be made, two distinct names, corresponding  to the two ways of determining the content, must be assigned to what these ways determine. 
The judgment, however, requires for 
its expression a sign for identity of content, a sign that connects these two names.~(...)

\quad Now let 
$$\Fa[1] (A \equiv B)$$   
mean that \textit{the sign} \textit{A} and \textit{the sign} \textit{B} \textit{have the same conceptual content}, \textit{so that we can everywhere put \textit{B} for \textit{A} and conversely}  \citep[\S8]{begr}. 
\end{quote}
\mh From the viewpoint of the distinction between sense and reference established in \textit{SR} (1892),  
the two names of the same point  mentioned  in the quotation above  have different senses but the same reference. Based on the \S8 of \textit{BS} we can say that 

\begin{description}
\item[] (4) two singular terms \textit{a} and \textit{b} have the same content if and only if  they 
pick up the same object. 
\end{description}
\mh    
Indeed, the identity of content in the proposition 52 of \textit{BS}  \citep[\S20]{begr}, 
\begin{description}
\item[] $ c \equiv d \to (f(c) \to f(d)) $,
\end{description}
\mh 
when \textit{c} and \textit{d} are singular terms, works extensionally as an identity axiom. 
It is to be noted that in the proposition 52, \textit{c} and \textit{d} can also be sentences, 
and $f$  any context, including an empty context, and so  
\begin{description}
\item[]  $ c \equiv d \to (c \to d) $
\end{description}
\mh can be derived. 
But in this case, \textit{c} and \textit{d} being sentences, 
Duarte has shown in \citep[pp. 334ff]{duarte} that $ c \equiv d $ 
does not follow from $c \to d $ and $d \to c $. 
Since the implication of \textit{BS} is the material implication,   
$ \equiv $ cannot be an extensional operator when flanked by sentences.

\subsection{Definitions in \textit{BS}}

The sign of identity of content also appears in \S24 of \textit{BS}, where Frege introduces his notation for definitions: $ {\Vdash } \ (A \equiv B)$. 
\begin{quote}
[A definition] differs from the judgments considered up to now in that it contains signs that have not been defined before; it itself gives the definition. It does not say ``The right side of the equation has the same content as the left'', but ``It is to have the same content'' \citep[\S24]{begr}.
\end{quote}
\mh 
A definition stipulates that expressions \textit{A} and \textit{B} have the same conceptual content.

\subsection{The problem with conceptual contents in \textit{BS}}	

We have just seen that the notion of identity of content is explained intensionally in \S3 
and extensionally in \S8. 
The symbol $ \equiv $ appears  in \S8 as an extensional operator, but 
later, in  \S24,  it is  an intensional operator. In  proposition 52, \S20, the symbol 
$ \equiv $ is extensional w.r.t. singular terms, and intensional w.r.t. sentences. 
Now, let us consider the sentences 
  
\begin{description}
\item[]	(5) The author of Nicomachean Ethics is Greek, 
\end{description}
\mh 
and 
	\begin{description}
	\item[]	(6) The tutor of Alexander the Great is Greek.
	\end{description}
\mh
Since the expressions `The author of Nicomachean Ethics' and `The tutor of Alexander the Great' pick up the same individual, Aristotle, according to (4), the sentences (5) and (6) have the same conceptual content. However, according to (3), they don't. 
Clearly, (5) and (6) are not intersubstitutable preserving correctness of inference because their consequences are not the same. 
The notion of conceptual content presented in 1879, in sections \S3 and \S8 of \textit{BS}, 
are incompatible because between sentences $ \equiv $ is an intensional operator, 
while flanked by singular terms it is an extensional operator.
Therefore, conceptual content cannot be the semantic value of the expressions of Frege's formal language.

\subsection{Frege's slingshot} \label{sling}

The name `slingshot' was given by Barwise and Perry  to a family of short arguments based on a small number of principles that intend to undermine important philosophical theses with just a few 
assumptions  \citep[cf.][pp. 375-378]{barper}. 
Different versions of the slingshot were used for different purposes, such as defending Frege's thesis that the reference of a sentence is its truth-value (cf. Church in  \cite[pp. 25-26]{church}), attacking  theories of facts (cf. Davidson in \cite[p. 42]{davidson}), and defending extensionality in general  (cf. Quine in \cite[p. 159]{quine}).

Neale in \cite{neale}  
works out a detailed analysis of several versions of the slingshot as formal arguments that prove that if a sentential connective $ \phi $ accepts in its scope certain principles of inference, $ \phi $ ends up being extensional. 
If $ \phi $ is presumably intensional, like the connectives that operate on facts, propositions, necessity, etc., according to the arguments presented by  Quine, Davidson and Church, 
if in the scope of $ \phi $ coextensional singular terms and logically equivalent sentences are intersubstitutable and the respective inferences are truth-preserving, it can be proved that $ \phi $ is extensional.  
Once the assumptions are accepted, applied to a theory of facts, the argument  proves that all facts collapse into a single fact \citep{davidson}; 
applied to a theory of propositions, the argument proves that all true propositions are synonymous -- or at least, as Church argues in  \citep{church}, all true propositions have the same reference.

\

The Fregean notion of conceptual content yields a slingshot argument that collapses all true identity sentences of arithmetic. In addition to (4), we need only to accept that

\begin{description}
\item[] 	(7) `\textit{b} is the successor of \textit{a}' and `\textit{a} is the predecessor of \textit{b}' have the same conceptual content, 
\end{description} 
\mh 
which seems to be perfectly justified by the example given by Frege in \S3 of \textit{BS}, quoted above, 
and also by the fact that the operations \textit{successor} and \textit{predecessor} are interdefinable, and so, according to \S24 of \textit{BS}, they have the same content. 
Now, given (4) and (7), all the sentences below have the same conceptual content: 

\begin{enumerate}  

\item[]	[1] 0 = the predecessor of 1 
\item[] [2] 1 =  the successor of 0  \quad\quad\quad	from [1] and (7)  
\item[] [3] 1 = the predecessor of 2 	\quad\quad from [2] and (4)
\item[]	[4] 2 = the successor of 1  \quad\quad\quad 		from [3] and (7)
\item[] and so on. 
\end{enumerate}

\mh 
Further substitutions can be made, and so we get that all true identity sentences $ a = b $ 
of Arithmetic have the same conceptual content.
Unlike other versions of the slingshot, this argument does not result in the collapse of all conceptual contents, but it is enough to reject the notion of conceptual content so defined as the semantic value of the expressions of Frege’s formal language.
Although 
Frege has not explicitly acknowledged this  argument,
it obviously gives support to
the claim that the reference of a sentence is its truth-value, since what is common to all these sentences is precisely the fact that they are true.

\section{A reconstruction of Frege's argument} \label{recon}

In this section, the argument in defense of truth-values as references will be analyzed and reconstructed. 
But let me begin by taking a look at the thesis that concepts are functions whose values are truth-values, 
presented in 1891 in the text \textit{Function and Concept} \citep{FC}. There, Frege extends the notion of 
function to include  expressions  
formed with the symbols $ = $, $ > $  and $ < $ 
\citep[p. 30]{FC}. 
When Frege asks what would be the values of such functions, the thesis that the reference of a sentence is a truth-value appears in the form of the claim that concepts are functions whose values are truth-values. 
Note that these two claims are virtually the same. The semantic value (i.e. reference) of the expression $ x^2 = 1 $  is a function. The value of that function, say, for the argument 2, will be the semantic value of the expression $ 2 \cdot 2 = 1 $, the truth-value \textit{false}. The claim that the values of a function like $ x^2 = 1 $, i.e. a concept, for different arguments are truth-values and the claim that the references of sentences  are truth-values are one and the same.

\

	The thesis that the reference of a sentence is its truth-value is presented and defended in \textit{SR} \citep[pp. 62-64]{SR}. 
Frege also argues in defence of this thesis in the posthumous 
\textit{Logic in Mathematics} \cite[pp. 231-233]{LM}  and in the letter to Russell 
\citep[pp. 152-153]{letter}. 
The most important difference between these arguments is the role of the principle of intersubstitutivity
of expressions with the same reference:

\begin{itemize} 
\item[] 	(IR) If \textit{Ref(A) = Ref(B)}, then \textit{Ref(...A...) = Ref(...B...)}\footnote{Note that (IR) is a consequence of the principle of compositionality. Furthermore, an analogous principle also holds for senses, since the sense of the whole expression depends on its structure and on the senses of its parts. In other words, both intensions and extensions are expected to behave in a compositional manner. But the principle of compositionality \textit{does not hold} for a `hybrid' notion, like the conceptual contents -- this is precisely the lesson from the slingshot argument.}.
\end{itemize}

\mh
In \textit{SR}, (IR) is a test applied by Frege after reaching the conclusion that the reference of a sentence is its truth-value. Later, in  1902 \citep{letter} and 1914 \citep{LM}, it is clear that (IR) is a premise, essential to justify the last step of the argument. 
	
	\medskip
	
Frege's argument  proceeds through three main steps: 
	
\begin{description}
\item[]	(I) Frege rules out the thought as the reference of sentences; 
\item[] (II) Frege concludes that sentences have reference; 
\item[] (III) Frege concludes that the reference of a sentence is its truth value. 
\end{description}	
\mh The critical step is the third. Let us take a closer look at these three steps.

\subsection{Step I}  
In \textit{SR}, after establishing the distinction between sense and reference for  names, Frege asks if such a distinction should be extended to complete sentences. 
Frege assumes that a sentence expresses a thought and asks if that thought should be considered 
to be the reference of the sentence. 
The answer is negative. Thoughts as references of sentences are rejected based on (IR) and a counterexample \citep[p. 62]{SR}. Frege considers, as he had already done earlier in the same text, that the sentences below

\begin{description}
\item[] 
(8) The morning star is a body illuminated by the Sun,
\end{description}
\mh and 
\begin{description}
\item[] 
(9) The evening star is a body illuminated by the Sun, 
\end{description}
\mh 
have different cognitive values, and so express different thoughts. But if the reference of a sentence were the thought expressed, the sentences (8) and (9)  should express the same thought, since the expressions `the morning star' and `the evening star' have the same reference. 
Hence, the thought cannot be the reference of a sentence.
 
\subsection{Step II}  
The second step aims to answer whether or not sentences have reference. The answer will be affirmative 
\citep[pp. 62-63]{SR}. 
It is worth noting that from this part of the text of \textit{SR}, together with \citep[pp. 152-153]{letter} and \citep[pp. 231-233]{LM}, the following equivalences may be established: 

\begin{description}
\item[] The parts of a sentence have reference \textit{if and only if }
\item[]  the complete sentence has reference \textit{if and only if }
\item[]  the sentence has a truth value \textit{if and only if }
\item[]  the thought expressed is true or false. 
\end{description}

\mh So, since we are interested in sentences with a truth-value, it follows that sentences have reference.
This conclusion, however, can be obtained from a short argument, as follows: 

\begin{description}
\item[] (10) If the parts of a sentence have reference, the complete sentence has reference. 
\item[] (11) If we are interested in truth, (we require that) the parts of a sentence have reference.
\item[] So, since in a scientific investigation we are interested in truth,   
\item[] (12) the complete sentence has reference. 
\end{description}
\mh 
The premises (10) and (11) are justified by the following passages: 

\begin{quote}
The fact that we concern ourselves at all about the reference of a part of the sentence indicates that we generally recognize and expect a reference for the sentence itself. The thought loses value for us as soon as we recognize that the reference of one of its parts is missing \citep[p. 63]{SR}.  
\end{quote}

\begin{quote}
What we talk about is the reference of words. We say something about the reference of the word `Sirius' when we say: `Sirius is bigger than the sun'. 
This is why in science it is of value to us to know that the words used have a reference. 
(...) The question first acquires an interest for us when we take a scientific attitude (...) Now it would be impossible to see why it was of value to us to know whether or not a word had a reference if the whole proposition did not have a reference and if this reference was of no value to us  \citep[p. 152]{letter}.
\end{quote} 

\begin{quote}
[it is essential] that the name `Etna' should have a reference, for otherwise we should be lost in fiction. The latter of course is essential only if we wish to operate in the realm of science. (...) 

If therefore we are concerned that the name `Etna' should designate something, we shall also be concerned with the reference of the sentence as a whole. That the name should designate something matters to us if and only if we are concerned with truth in the scientific sense. \citep[p. 232]{LM}.
\end{quote}

\mh At this point, it is important to call attention to the fact that Frege's conclusion \textit{is not}
that sentences designate something in the same sense that names designate objects. 
Frege must be understood as saying that, given that the parts of a sentence 
of his formal language have a semantic value, the complete sentence also has a semantic value.
A possible objection here is that since 
Frege did not have a semantic theory, as we understand it today, 
Frege could not be asking for semantic values.
It is true that semantics became a systematic  
discipline some decades later, mainly in the work of Tarski on model theory.
Nevertheless, Frege's formal language was an \emph{interpreted} language. 
Recall that  Frege opposed himself to formalism in Mathematics, and so, 
Arithmetic could not be a manipulation of empty or meaningless symbols. 
But this is virtually the same as saying that he had to attribute semantic values 
to the expressions of his formal language, conceived precisely to deal with  Arithmetic. 
Besides, as far as Frege was talking about reference and truth -- not by chance the 
basic notions of the (extensional) semantics of classical logic -- he was using semantic notions in our sense of semantics.

\subsection{Step III} 
After concluding that sentences have reference, in \textit{SR}, Frege says: 

\begin{quote}
We have seen that the reference of a sentence may always be sought, whenever the reference of its components is involved; and that this is the case when and only when we are inquiring after the truth value \citep[p. 63]{SR}.
\end{quote}
\mh 
And concludes, somewhat reluctantly, that

\begin{quote}
We are therefore driven into accepting the \textit{truth value} of a sentence as constituting its reference. By the truth value of a sentence I understand the circumstance that it is true or false (Ibidem).
\end{quote}

\mh Two paragraphs later \citep[p. 64]{SR} Frege applies (IR) as a test to confirm, 
or to make more plausible, his conclusion: 

\begin{quote}
If our supposition that the reference of a sentence is its truth value is correct, the latter 
must remain unchanged when a part of the sentence is replaced by an expression 
having the same reference. And this is in fact the case \citep[p. 64]{SR}. 
\end{quote}

\mh However, in both  \cite{letter} and \citep{LM}, (IR) is a premise of the argument. 

\begin{quote}
The reference of the proposition [sentence] must be something which does not change when one sign is replaced by another with the same meaning but a different sense. What does not change in the process is the truth-value  \citep[p. 152]{letter}.
\end{quote}

\begin{quote}
The reference of a sentence must be something which remains the same, if one of the parts is replaced by something having the same reference. We return now to the sentence `($16 - 2$) is a multiple of 7'. (...) 

[W]hat is not altered by replacing the sign `($16 - 2$)' by the sign `($17 - 3$)', whose reference 
is  the same, is what I call the \textit{truth-value}.~(...)

We say accordingly that sentences have \textit{the same reference} if they are both true,
or if they are both false.  \citep[p. 232-233]{LM}.
\end{quote}

In my view, the best way to read Frege's argument is as follows. Frege concluded that sentences have reference, and that the truth-value of a sentence is a plausible candidate, both because it shows itself as an alternative in the argument that concludes that sentences have reference, and also because it satisfies (IR), a necessary condition for being the reference. From this, Frege draws the conclusion that the reference of a sentence is its truth-value. So, the argument is

\begin{description}
\item[] 
(13) The truth-value is a plausible candidate for being the reference of a sentence;
\item[]
(14) If something is the reference of a sentence, it must satisfy (IR);
\item[]
(15) The truth-value satisfies (IR);
\item[] 
Therefore:
\item[] 
(16) The reference of a sentence is its truth-value.
\end{description}

\mh The problem is that the above argument has not explicitly excluded something different from truth-values that could play the role of reference. Indeed, the step from (14) and (15) to (16) is invalid -- it is an instance of the fallacy of affirmation of the consequent. Although Frege does not explicitly claim that there is no third option besides the truth-value and the thought for being the reference of a sentence, it is clear from the quotations from \cite{letter} and \citep{LM}  that he was convinced that there was no third alternative.

\subsection{The reconstructed argument} 

If we put the whole scenario into perspective, from 1879 (\textit{BS}) to 1892 (\textit{SR}), 
we see three options for being the semantic value of a sentence of Frege's formal system: 
(i) the conceptual content,  (ii) the thought expressed, or  (iii) the truth-value.
The conceptual content is excluded because of the problems discussed in  section \ref{conc.cont} above. Thus, the argument can be reconstructed by adding  premises (17) and (18) below:  

\begin{description}
\item[] (14) If something is the reference of a sentence, it must satisfy (IR);
\item[] (15) The truth-value satisfies (IR);
\item[] (17) The thought does not satisfy (IR); 
\item[] (18) The reference of a sentence is either the truth-value or the thought;
\item[] Therefore:
\item[] (16) The reference of a sentence is its truth-value. 
\end{description}

Let us consider now a possible objection, namely, that there could still be a notion suited to be 
the reference but different from truth-values, thoughts, and conceptual contents. 
Let's call such a notion $ \Theta $, and it would be something inbetween the thought and the truth-value, 
less intensional than the thought, but not extensional like the truth-value. 
In this case, $ \Theta $ would still imply that the sentences `\textit{b} is the successor of \textit{a}' 
and `\textit{a} is the predecessor of \textit{b}' have the same reference, 
since the first can be defined from 
the second (or vice-versa), and it is not plausible to say that in a definition 
the definiendum and the definiens do not have the same reference. 
So,  again, all true identity sentences of Arithmetic would have the same reference. 
But a context where any true atomic sentence may be substituted by any other true atomic sentence 
is nothing but an extensional context, that is, 
a context in which the semantic value of a sentence is its truth-value.

\section{Final remarks}

The thesis that truth-values are references of sentences became the standard approach in the semantics of first-order logic, as we learn from any book of elementary logic. 
My attempt to reconstruct Frege's argument in defence of this thesis and to answer the criticisms 
made by Russell, Chateaubriand and Dummett \citep[cf.][]{russell, dummett73, chat} has been 
motivated by the feeling that these  views, besides being mistaken,  conceal the  two important achievements of Frege's works on logic mentioned in the Introduction: (i) that an extensional logic is very suitable as an account of logical
consequence for Mathematics, and (ii) that  in such an extensional logic, the
extension of a sentence is its truth-value.

My strategy here has been to investigate the problem from a viewpoint that considers the development of 
Frege's doctrines from 1879 to 1893. I have  tried to show that the argument that concludes that 
the reference of a sentence is its truth-value is not really in \textit{SR}, 
but rather in the path that goes from \textit{BS} to \textit{BLA}. 
And it depends on the distinction established by Frege in \textit{SR} 
between the intension and the extension of the expressions of a formal language.  

We have seen that a problem of Frege's argument in defence of truth-values as references is that 
it seems that he has not considered a possible third alternative besides 
thoughts and  truth-values.
Actually,  
what shows that Frege did not have a third alternative is the slingshot argument that collapses conceptual contents of true identities (Section \ref{sling} above), together with the problems of the notion of conceptual content. 
The latter is simultaneously intensional and extensional, respectively, 
to sentences and to singular terms. Note that such a mixture of intensionality with extensionality is 
directly responsible for the collapse that results from the slingshot argument. 
   
The thesis that truth-values are references of sentences is indeed not plausible 
from the viewpoint of a theory of meaning. 
But this  has nothing to do with the purpose of providing an account of logical consequence 
and a formal language to be used 
in the logicist project. 
A theory of meaning and a theory of logical consequence are two quite different things. 
For the reader not familiar with the role of \textit{SR} within Frege's work, it may indeed seem that 
the discussion about sense and reference, and the detailed analysis of intensional contexts 
in the second part of \textit{SR},  
is a contribution for a theory of meaning -- and in fact it is, 
but this is a secondary issue, not the primary interest of Frege. 
The supposed implausibility of truth-values as references, alleged by Chateaubriand \citep{chat}, Dummett 
\citep{dummett73}, and Russell \citep{russell}, disappears when we realize that the central point of \textit{SR}, to present and defend the thesis that the reference of a sentence is its truth value,  
is only  to establish that the logic of \textit{BLA}, 
to be published  the next year, is extensional.
Indeed,  what is designated by a singular term coincides with its semantic value.
Maybe  this fact has an important role in the misreadings of Frege mentioned here, 
but of course, Frege was not  interested in whether or not a sentence designates a fact, a state of affairs, or anything like that.\footnote{The author would like to thank Alessandro Duarte, Dirk Greimann,  the audience of the \emph{Sixth International Symposium on Philosophy of Language and Metaphysics} (Universidade Federal Fluminense, Rio de Janeiro,  November  2017), and two  anonymous  referees  for some valuable comments on a previous version of this text.}



%

\end{document}